\documentclass[11pt]{amsart}

\renewcommand{\bar}{\overline}
\newcommand{\eps}{\epsilon}
\newcommand{\pa}{\partial}
\newcommand{\ph}{\varphi}

\font\strange=msbm10

\newcommand{\C}{{{\mathchoice  {\hbox{$\textstyle{\text{\strange C}}$}}
{\hbox{$\textstyle{\text{\strange C}}$}}
{\hbox{$\scriptstyle{\text{\strange C}}$}}
{\hbox{$\scriptscriptstyle{\text{\strange C}}$}}}}}

\newcommand{\R}{{{\mathchoice  {\hbox{$\textstyle{\text{\strange R}}$}}
{\hbox{$\textstyle{\text{\strange R}}$}}
{\hbox{$\scriptstyle{\text{\strange R}}$}}
{\hbox{$\scriptscriptstyle{\text{\strange R}}$}}}}}

\newcommand{\Z}{{{\mathchoice  {\hbox{$\textstyle{\text{\strange Z}}$}}
{\hbox{$\textstyle{\text{\strange Z}}$}}
{\hbox{$\scriptstyle  Z\kern-0.3em  Z$}}
{\hbox{$\scriptscriptstyle  Z\kern-0.2em  Z$}}}}}

\newcommand{\N}{{{\mathchoice  {\hbox{$\textstyle{\text{\strange N}}$}}
{\hbox{$\textstyle{\text{\strange N}}$}}
{\hbox{$\scriptstyle  N\kern-0.3em  N$}}
{\hbox{$\scriptscriptstyle  N\kern-0.2em  N$}}}}}

\newcommand{\frk}[1]{{\mathfrak{#1}}}

\renewcommand{\bar}{\overline}

\usepackage{amsmath,amsthm,amscd}

\title
[]{On the Rigidity of  Horizontal Slices}
\author[]{
Zhiqin Lu}
\date{\today}

\address[Zhiqin Lu]
{Department of Mathematics\\
University of California, Irvine\\
Irvine, CA 92697}
\email{zlu@math.uci.edu}
\thanks{Research supported by  NSF grant DMS 0196086.}

\newtheorem*{thm}{Theorem}
\newtheorem{theorem}{Theorem}[section]
\newtheorem{lemma}{Lemma}[section]
\newtheorem{cor}{Corollary}[section]
\newtheorem{prop}{Proposition}[section]

\newtheorem{definition}{Definition}[section]

\newtheorem{assumption}{Assumption}[section]

\theoremstyle{remark}
\newtheorem{rem}{Remark}[section]

\begin{document}
\begin{abstract}
In this paper, we proved a rigidity theorem
of the Hodge metric for concave horizontal slices and
a local rigidity theorem for the monodromy representation.
\end{abstract}

\maketitle

\numberwithin{equation}{section}


\section{Introductions}
Let $(X,\omega)$ be a polarized simply connected Calabi-Yau manifold.
That is, 
$X$ is an $n$-dimensional compact K\"ahler 
manifold with zero first Chern class
and $[\omega]\in H^2(X,\Z)$ is a K\"ahler metric. By 
the famous theorem of Yau~\cite{Y1}, there is 
a K\"ahler metric on $X$ 
in the same cohomological class of $[\omega]$
such that its Ricci 
curvature  is zero.

Let $\Theta$ be the holomorphic tangent 
bundle of $X$. In ~\cite{T1}, Tian proved that
the universal deformation space of the complex structure is smooth. 
The complex
dimension of the universal deformation space is 
$dim\,H^1(X,\Theta)$.
 In other words, there are no 
obstructions towards the deformation 
of the complex structure of  Calabi-Yau manifold. A good
reference for the proof is in ~\cite{Fd}.

Take $n=3$ for example.
 A natural 
question is that to what extent the Hodge 
structure, namely, the decomposition of 
$H^3(X,\C)$, into the sum of $H^{p,q}$'s 
(p+q=3), determines a Calabi-Yau threefold. 
Let's recall 
 the concept of classifying space in~\cite{Gr}, 
which is a generalization of classical 
period domain. In the case of Calabi-Yau 
threefold, the classifying space $D$ is 
defined as the set of the filtrations of 
$H=H^3(X,\C)$ by
\[
0\subset F^3\subset F^2\subset F^1\subset H
\]
with $dim\,F^3=1$, 
$dim\,F^2=n=dim\,H^1(X,\Theta)$, 
$dim\,F^1=2n+1$, and $H^{p,q}=F^p\cap 
\bar{F}^q, H=F^p\oplus \bar{F^{4-p}}$ (p+q=3) together with a quadratic 
form $Q$ such that 
\begin{enumerate}
\item $i\, Q(x,\bar{x})<0$
if $0\neq x\in H^{3,0}$
\item $i\, Q(x,\bar{x})>0$
if $0\neq x\in H^{2,1}$
\end{enumerate}
where $i=\sqrt{-1}$.

There is a natural map from the universal deformation space into the 
classifying space. Intuitively, this is because $D$ is 
just the set of all the possible ``Hodge 
decompositions". 
Such a map is called a period map.
In the case of Calabi-Yau, 
the map is a holomorphic immersion. Thus  in 
that case, the infinitesimal Torelli theorem 
is valid~\cite{Gr}.

It can be seen that $D$ fibers over a 
symmetric space $D_1$. But such a symmetric 
space needs not to be Hermitian. Even $D_1$ 
is Hermitian symmetric, $D$ still needs not 
fiber holomorphically over $D_1$. Although in 
that case, there is a complex structure on $D$ 
such that $D$ becomes homogeneous 
K\"ahlerian~\cite{M}. 

Griffiths introduced the concept of 
horizontal distribution in~\cite{Gr}.
He proved that the image of the universal deformation space, via the
period map to the
classifying space, is an integral submanifold of the horizontal
distribution. A horizontal slice is an integral 
complex submanifold of the
horizontal distribution. In this terminology, the 
universal deformation space is a
horizontal slice of
the classifying space. The horizontal 
distribution is a highly nonintegrable system.

Because of the above result of Griffiths,
it is   interesting to study 
horizontal slices of a classifying space. The local properties of
horizontal
slices have been studied in~\cite{Lu6} and~\cite{Lu5}.

In~\cite{Lu6},
we introduced  a new K\"ahler metric on a horizontal slice
$U$. We call such a metric the Hodge metric. The main result in~\cite{Lu6}
is that (see also~\cite{Lu5} for the case $n=3$)

\begin{thm}
Let $U\rightarrow D$ be a horizontal slice. Then the restriction
of the natural invariant Hermitian metric of $D$ to $U$ is actually
K\"ahlerian. We call such a K\"ahler metric the Hodge metric of $U$. The
holomorphic bisectional curvature of the Hodge metric is nonpositive. The
Ricci curvature of the Hodge metric is negative away from zero.
\end{thm}

In this paper, we study some global rigidity properties of horizontal
slices. 
In order to do that, we observe that the universal deformation space $U$
carries less global information than the moduli space ${\mathcal M}$,
which is essentially the quotient of the universal deformation space by a
discrete subgroup $\Gamma$ of $Aut(U)$. The group $\Gamma$ is called the
monodromy group of the moduli space.
The volume of the space $\Gamma\backslash U$ is finite with respect to the
Hodge metric. This is the consequence of the theorem of Viehweg~\cite{V2},
the theorem of Tian~\cite{T7} and the above theorem.

For a horizontal slice $U$ of $D$, if $\Gamma$ is a discrete subgroup of
$U$ such that the volume of $\Gamma\backslash U$ is finite, then a general
conjecture is that whether $\Gamma$ completely determines the space
$\Gamma\backslash U$. In the case where $n=2$, this is correct by the
superrigidity theorem of Margulis~\cite{Mar}. In general, this is a very
difficult
problem.

We consider the following weaker rigidity  problem: to what extent 
the complex
structure of the moduli space determines the metrics on the moduli space
and the monodromy representation?  To this problem, we have the 
following result in this paper.

First, we proved that, if $\Gamma\backslash U$ is a complete
concave manifold, then the complex structure of $\Gamma\backslash U$
completely determines the Hodge metric of $\Gamma\backslash U$. More
precisely,
 we proved that if for some discrete group $\Gamma$ of
$Aut(U)$, $\Gamma\backslash U$
is a concave complete complex manifold, then the Hodge metric defined
on
$U$ is intrinsic. In other words, the Hodge metric doesn't depend on the
choice of the holomorphic immersion $U\rightarrow D$ 
from which it becomes 
a horizontal slice.

\begin{theorem}
If the moduli space $\Gamma\backslash{U}$ is a  concave 
manifold. Then the Hodge metric is intrinsically defined. 
\end{theorem}

The second main result of this paper is the local rigidity of the
monodromy representation. The result is in  Theorem~\ref{thm51}. 
We combine the superrigidity theorem of Margulis~\cite{Mar} 
together with some ideas of Frankel~\cite{SF} in the proof of the theorem.

\begin{theorem}\label{thm51} (For definition of the notations,
see 
\S 5) Let $\Gamma\backslash{U}$ be of finite Hodge volume. 
Suppose further that ${\mathcal G}_0$
is semisimple 
and ${\mathcal G}_0/{\mathcal K}_0$ is a Hermitian symmetric space but is
not a complex ball, where ${\mathcal K}_0$ is the maximum compact subgroup
of ${\mathcal G}_0$.  Then the representation
 $\Gamma\rightarrow\,{\mathcal G}$ is 
locally rigid.
\end{theorem}

The motivation behind the above results is that in the case of 
$K$-$3$ surfaces, the moduli space is a local symmetric space of rank 2.
 But even in the case of Calabi-Yau
threefold, little has been known about the moduli space. We wish to
involve certain kinds of 
metrics (Weil-Petersson metric, Hodge metric, etc) in the study of the
moduli space of Calabi-Yau manifolds.
The metrics have applications in Mirror Symmetry of Calabi-Yau
manifolds~\cite{Y2}.

{\bf Acknowledgment.} The author thanks Professor G.  Tian for
his help and
encouragement during the preparation of this paper.

\section{Preliminaries}
In this section, we give some definitions and notations which will be 
used throughout this paper.
Unless otherwise stated, the materials in this section are from the book
of Griffiths~\cite{Gr}.

Let $X$ be a compact K\"ahler manifold. A $C^\infty$ form on $X$ 
decomposes into (p,q)-components according to the number of 
$dz's$ and $d\bar{z}'s$. Denoting the $C^\infty$ $n$-forms and
the 
$C^\infty(p,q)$ forms on $X$ by $A^n(X)$ and $A^{p,q}(X)$ 
respectively, we have the decomposition
\[
A^n(X)=\underset{p+q=n}{\oplus}A^{p,q}(X)
\]

The cohomology group is defined as
\begin{align*}
H^{p,q}(X)=&\{closed (p,q) -forms\}/\{exact (p,q)-forms\}\\
=&\{\phi\in A^{p,q}(X)|d\phi=0\}/dA^{n-1}(X)\cap A^{p,q}(X)
\end{align*}

The following theorem is well known:

\begin{thm}[Hodge Decomposition Theorem]   Let $X$ be a 
compact K\"ahler manifold of dimension $n$. Then the n-th 
complex de Rham cohomology group of $X$ can be written as a 
direct sum
\begin{equation*}
H^n_{DR}(X,\Z)\otimes {\C}
=H^n_{DR}(X,\C)=
\underset{p+q=n}{\oplus}H^{p,q}(X)
\end{equation*}
such that $H^{p,q}(X)=\bar{H^{q,p}(X)}$.
\end{thm}

\begin{rem}
We can define a filtration of $H^n_{DR}(X,\C)$ by
\[
0\subset F^n\subset F^{n-1}\subset \cdots F^1=H=H^n_{DR}(X,\C)
\]
such that
\[
H^{p,q}(X)=F^p\cap\bar{F}^q
\]
So the set $\{H^{p,q}(X)\}$ and $\{F^p\}$ are equivalent in defining 
the Hodge decomposition. In the remaining of this paper, we 
will use both notations interchangeably.
\end{rem}

\begin{definition}
A Hodge structure of weight $j$, denoted by $\{H_{Z},H^{p,q}\}$, is given
by
a
lattice $H_{Z}$ of finite rank together with a decomposition on its
complexification $H=H_{Z}\otimes \C$
\[
H={\displaystyle \underset{p+q=j}{\oplus}} H^{p,q}
\]
such that
\[
H^{p,q}=\bar{H^{q,p}}
\]
\end{definition}

A polarized algebraic manifold is a pair $(X,\omega)$ consisting
of an algebraic manifold $X$ together with a K\"ahler form $\omega\in
H^2(X,\Z)$. Let
\[
L: H^j(X,\C)\rightarrow H^{j+2}(X,\C)
\]
be the multiplication by $\omega$, we recall below two fundamental
theorems of Lefschetz:

\begin{thm}[Hard Lefschetz Theorem] On a polarized algebraic manifold
$(X,\omega)$ of dimension $n$,
\[
L^k : H^{n-k}(X,\C)\rightarrow H^{n+k}(X,\C)
\]
is an isomorphism for every positive integer $k\leq n$.
\end{thm}

From the theorem above, we know that
\[
L^{n-j}: H^j(X,\C)\rightarrow H^{2n-j}(X,\C)
\]
is an isomorphism for $j\geq 0$. The primitive cohomology $P^j(X,\C)$ is
defined to be
the kernel of $L^{n-j+1}$ on $H^j(X,\C)$.

\begin{thm}[Lefschetz Decomposition Theorem]
On a polarized algebraic manifold $(X,\omega)$, we have for any integer
$j$ the following decomposition
\[
H^j(X,\C)=\underset{k=0}{\overset{[\frac n2]}\oplus}
L^k P^{j-2k}(X,\C)
\]
\end{thm}

It follows that the primitive cohomology groups determine completely the
full complex cohomology.

In this paper  we are only interested in the cohomology group
$H^n_{DR}(X,\C)$. Define
\[
H_Z=P^n(X,\C)\cap H^n(X,\Z)
\]
and 
\[
H^{p,q}=P^n(X,\C)\cap H^{p,q}(X)
\]

Suppose that $Q$ is the quadric form on $H^n_{DR}(X,\C)$ 
induced by the cup product of the cohomology group. 
$Q$ can be represented by
\[
Q(\phi, \psi)=(-1)^{n(n-1)/2}\int \phi\wedge \psi
\]

$Q$ is a nondegenerated form, and is skewsymmetric if $n$ is odd 
and is symmetric if $n$ is even. It satisfies the two 
Hodge-Riemannian relations
\begin{enumerate}\label{hr}
\item $Q(H^{p,q},H^{p',q'})=0 \quad unless\quad p'=n-p, q'=n-q$;
\item $(\sqrt{-1})^{p-q}\,Q(\phi,\bar\phi)>0$ for any nonzero 
element $\phi\in H^{p,q}$.
\end{enumerate}

Let $H_Z$ be a fixed lattice, $n$  an integer, $Q$ a bilinear
form on $H_Z$, which is symmetric if $n$ is even and
skewsymmetric if $n$ is odd. And let  $\{h^{p,q}\}$ be a collection of
integers such that $p+q=n$ and $\sum h^{p,q}=rank\, H_Z$.
Let $H=H_Z\otimes \C$.

\begin{definition}
A polarized Hodge structure
of weight $n$, denoted by\newline
$\{H_Z,F^p,Q\}$, is given by a 
filtration of $H=H_Z\otimes \C$
\[
0\subset F^n\subset F^{n-1}\subset\cdots\subset F^0\subset H
\]
such that
\[
H=F^p\oplus\bar{F}^{n-p+1}
\]
together with a bilinear form
\[
Q: H_Z\otimes H_Z\rightarrow \Z
\]
which is skewsymmetric if $n$ is odd and symmetric if $n$ is 
even such that it satisfies the two Hodge-Riemannian 
relations:
\begin{enumerate}\label{hrf}
\item $Q(F^p,F^{n-p+1})=0\quad unless\quad p'=n-p,q'=n-q$;
\item $(\sqrt{-1})^{p-q}\,Q(\phi,\bar\phi)>0$ if $\phi\in 
H^{p,q}$ and $\phi\neq 0$
\end{enumerate}
where $H^{p,q}$ is defined by
\[
H^{p,q}=F^p\cap\bar{F}^q
\]
\end{definition}

\begin{definition}
With the notations as above, the classifying space $D$ for the 
polarized Hodge structure is the set of all the filtration
\[
0\subset F^n\subset\cdots\subset F^1\subset H, dim\, F^p=f^p
\]
with $f^p=h^{n,0}+\cdots +h^{n,n-p}$ on which $Q$ satisfies the 
Hodge-Riemannian relations as above.
\end{definition}

$D$ is a complex homogeneous space.
Moreover, $D$ can be written as $D=G/V$ where $G$ is a noncompact
semisimple Lie group and $V$ is its compact subgroup.
In general, $D$ is \textsl{not} a
homogeneous K\"ahler manifold.

\section{The Canonical Map and the Horizontal 
Distribution}
In this section we study some elementary properties of  classifying space
and  horizontal slice.

Suppose $D=G/V$ is a classifying space. We fix a point of $D$, say $p$, 
which can be represented by  the subvector spaces of $H$
\[
0\subset F^n\subset F^{n-1}\subset\cdots F^1\subset H
\]
or the set 
\[
\{H^{p,q}|p+q=n\}
\]
described in the previous section. 
We  define the subspaces of $H$:
\begin{align*}
& H^+=H^{n,0}+H^{n-2,2}+\cdots\\
& H^-=H^{n-1,1}+H^{n-3,3}+\cdots
\end{align*}

Suppose $K$ is the subgroup of $G$ such that $K$ leaves $H^+$ invariant. 
Then we have

\begin{lemma}
\label{l21}
The identity component $K_0$ of $K$ is the maximal connected compact
subgroup of 
$G$ containing $V$. 
In particular,  $V$ itself is a compact subgroup.
\end{lemma}

{\bf Proof:} Recall that $V\subset G\subset Hom(H_R,H_R)$ is a 
{\it real} subgroup, where $H_R=H_Z\otimes \R$.
Without losing generality, we assume
 $V$ fixes $p$. Then  we have
\[
VF^p\subset F^p
\]
where $p=1,\cdots,n$. This implies that
\[
V\bar{F^q}\subset\bar{F^q}
\]
for $q=1,\cdots,n$. So
\[
VH^{p,q}=V(F^p\cap\bar{F^q})\subset VF^p\cap V\bar{F^q}=H^{p,q}
\]

Thus $V$ leaves $H^+$ invariant and thus $V\subset K$.

In order to prove that $K_0$ is a compact subgroup,
we fix some $H^+,H^-\subset H$.
Note that if $0\neq x\in H^+$, then from the second Hodge-Riemannian
relation
\[
(\sqrt{-1})^n\,Q(x,\bar{x})>0
\]

So for any norm on $H^+$, there is a $c>0$ such that 
\[
\frac{1}{c}||x||^2\ge(\sqrt{-1})^n\,Q(x,\bar{x})\ge c||x||^2
\]

For the same reason, we have
\[
\frac{1}{c}||x||^2\ge -(\sqrt{-1})^n\,Q(x,\bar{x})\ge c||x||^2
\]
for $x\in H^-$.

Let $g\in K_0$. For any $x$, let $x=x^++x^-$ be the 
decomposition of $x$ into $H^+$ and $H^-$ parts. Then
\begin{align*}
&||gx^{\pm}||^2\le\pm\frac{1}{c}(\sqrt{-1})^n\,Q(gx^{\pm}
,g\bar{x^{\pm}})\\
&=\pm\frac{1}{c}(\sqrt{-1})^n\,Q(x^{\pm},\bar{x^{\pm}})\le\frac{1}{c^2} 
||x^{\pm}||^2 
\end{align*}

Thus
\[
||g||\le C
\]

So the norm of the element of $K_0$ is uniformly bounded. 
Consequently,
 $K_0$ is a compact subgroup.

Suppose that $K'\supset K_0$ is a compact connected subgroup. 
Suppose ${\frk k}'$ is the Lie algebra of $K'$, then if $K_0$
is not maximal, there is a $\xi\in 
{\frk k}'$ such that $\xi\notin {\frk f}_0$ for the Lie algebra 
${\frk f}_0$ of $K_0$. 

Suppose $\xi=\xi_1+\xi_2$ is the decomposition for which
\begin{align*}
& \xi_1 : H^+\rightarrow H^+, H^-\rightarrow H^-\\
& \xi_2 : H^+\rightarrow H^-, H^-\rightarrow H^+
\end{align*}

Then we have

\begin{lemma}
\label{lem22}
$\xi_1,\xi_2\in {\frk g}_R$ for the Lie algebra ${\frk g}_R$ of $G$.
\end{lemma}

 {\bf Proof:}  First we observe that
\begin{align*}
& Q(H^+,H^+)=Q(H^-,H^-)=0, & n \quad odd\\
&Q(H^+,H^-)=Q(H^-,H^+)=0, & n \quad  even
\end{align*}
by the type consideration. Since $Q$ is invariant under the 
action of $G$ by definition, we have
\[
Q(\xi x,y)+Q(x,\xi y)=0
\]
Thus
\[
Q(\xi_1 x,y)+Q(x,\xi_1 y)+Q(\xi_2 x,y)+Q(x,\xi_2 y)=0
\]
If $n$ is odd then if $x\in H^+,y\in H^+$ or $x\in H^-,y\in H^-$
then
\[
Q(\xi_1 x,y)+Q(x,\xi_1 y)=0
\]
so in this case
\[
Q(\xi_2 x,y)+Q(x,\xi_2 y)=0
\]
and if $x\in H^+,y\in H^-$ or $x\in H^-,y\in H^+$ then we have
\[
Q(\xi_2 x,y)+Q(x,\xi_2 y)=0
\]
automatically. Thus we concluded 
\[
Q(\xi_2 x,y)+Q(x,\xi_2 y)=0
\]
for any $x,y\in H$. So $\xi_2\in {\frk g}_R$ and thus $\xi_1\in{\frk g}_R$.

The same is true if $n$ is even. 

\qed

We define the Weil operator
\[
C: H^{p,q}\rightarrow H^{p,q},\quad C|_{H^{p,q}}=(\sqrt{-1})^{p-q}
\]
Then we have
\[
C|_{H^+}=(\sqrt{-1})^n,\qquad C|_{H^-}=-(\sqrt{-1})^n
\]
 
Let
\[
Q_1(x,{y})=Q(Cx,\bar{y})
\]
 
Then we have

\begin{lemma}
$Q_1$ is an Hermitian inner product.
\end{lemma}
 
{\bf Proof:}
Let
\[
x=x_1+x_2
\]
be the decomposition of $x$ such that $x_1\in H^+$ and $x_2\in
H^-$.
 
If $n$ is odd, then $\bar{x_2}\in H^+$. So
$Q(x_1,\bar{x_2})=0$; if
$n$ is even, then $\bar{x_2}\in H^-$. So $Q(x_1,\bar{x_2})=0$.

If $x\neq 0$ we have
\begin{align*}
Q_1(x,{x}) & =Q_1(x_1,{x_1})+Q_1(x_2,{x_2})
+Q_1(x_1,{x_2})+Q_1(x_2,{x_1})\\
&= Q_1(x_1,{x_1})+Q_1(x_2,{x_2})\\
&= (\sqrt{-1})^n(Q(x_1,\bar{x_1})-Q(x_2,\bar{x_2}))>0
\end{align*}
 
Thus $Q_1(\cdot,\cdot)$ is a Hermitian product on $H$. 
Furthermore, it defines an inner product on $H_R=H_Z\otimes \R$.

\qed

Now back to the proof of Lemma ~\ref{l21}, we have
\[
Q_1(\xi_2 x,y)=Q(C\xi_2 x,\bar{y})=-Q(\xi_2 Cx,\bar{y})=Q(Cx,\xi_2 \bar{y})
=Q_1(x,\xi_2 y)
\]

Thus $\xi_2$ is a Hermitian metrics under the metric $Q_1$.  
Since $K'$ is a compact group, there is a constant $C$ such that
\[
||exp\,(t\xi_2)||\le C<+\infty
\]
for all $t\in \R$ which implies $\xi_2=0$.

\qed

\begin{lemma}\label{lem34}
Let
\[
D_1=\{H^{n,0}+H^{n-2,2}+\cdots|\{H^{p,q}\}\in D\}
\]
Then the group $G$ acts on $D_1$ transitively with the stable subgroup 
$K_0$, and  $D_1$ is a symmetric space.
\end{lemma}

{\bf Proof:}
For $x,y\in D_1$, let $H_x^{p,q},H^{p,q}_y$ be the corresponding 
points in $D$. Since $D$ is homogeneous, we have a $g\in G$ such that
\[
g\{H_x^{p,q}\}=H_y^{p,q}
\]
So $gx=y$.
This proves that $G$ acts on $D_1$ transitively.
 By definition, $K_0$ fixes the $H^+$ of 
the fixed point $p\in D$. By Lemma~\ref{l21}, $D_1$ is a 
symmetric space.

\begin{definition}\label{pro}
We call the map $p$
\[
p: G/V\rightarrow G/K_0,\qquad \{H^{p,q}\}\mapsto H^{n,0}+H^{n-2,2}+\cdots
\]
the natural projection of the classifying space.
\end{definition}
\label{page8}
There
are universal holomorphic bundles $\underline{F}^n,
\cdots,\underline{F}^1, \underline{H}$ over $D$, namely we
assign any
point $p$ of $D$ the linear space
\[
0\subset F^n\subset\cdots\subset F^1\subset H
\]
or in other words, assign every point of $D$ the space
$H=H_{Z}\otimes \C$, with the Hodge decomposition
\[
H=\sum H^{p,q}
\]
 
It is well known that the holomorphic tangent bundle $T(D)$ can
be realized by
\[
T(D)\subset \oplus Hom(F^p,H/F^p)=
\underset{r> 0}{\oplus}
Hom(H^{p,q},H^{p-r,q+r})
\]
such that the following compatible condition holds
 
\[
\begin{CD}
F^{p} @> >>F^{p-1}\\
@VVV        @VVV\\
H/F^p@ <<< H/F^{p-1}
\end{CD}
\]
 
We  define a subbundle $T_h(D)$  called the
horizontal bundle of $D$, by
\[
T_h(D)=\{\xi\in T(D)| \xi F^p\subset F^{p-1}\}
\]
 
$T_h(D)$ is called the horizontal distribution of $D$.
 The properties of the horizontal bundle or the horizontal
distribution play an important role in the theory of moduli space.

Let ${\frk g}_R$ be the Lie algebra of $G$. Suppose
\[
{\frk g}_R={\frk f}_0+{\frk p}_0
\]
is the Cartan decomposition of ${\frk g}_R$ into the compact and
noncompact part. 

\begin{lemma}
If we identify $T_0(G)$ with the Lie algebra ${\frk g}_R$. Then
\[
E\subset {\frk p}_0
\]
where $E$ is the fiber of $T_h(D)$ at the original point.
\end{lemma}

{\bf Proof:}
Suppose 
\[
\{0\subset f^n\subset f^{n-1}\subset\cdots f^1\subset H\}\quad
or\quad\{h^{p,q}\}
\]
is the set of subspace representing the point $eV$
of $D=G/V$.
Suppose $X\in E$. Then $X\in E$ if
\[
X: f^k\rightarrow f^{k-1}
\]
Let $X=X_1+X_2$ be the Cartan decomposition
with $X_1\in{\frk f}_0$, and $X_2\in{\frk p}_0$. Let
\begin{align*}
& h^+=h^{n,0}+h^{n-2,2}+\cdots\\
& h^-=h^{n-1,1}+h^{n-3,3}+\cdots
\end{align*}
be the subspaces of $H$.

By definition $X_1\in{\frk f}_0$, we see that
\[
X_1: h^+\rightarrow h^+,\qquad h^-\rightarrow h^-
\]
Since $X$ maps $f^k$ to $f^{k-1}$,  so does $X_1$.
So $X_1$ must leave  $f^k$ invariant because $X_1$ sends $h^+$ to
$h^+$, and $h^-$ to $h^-$.

From the above argument we see that
 $X_1\in{\frk v}$, the Lie algebra of $V$. Thus 
the action $X$ 
 on the classifying space is the same as $X_2$.
But $X_2\in{\frk p}_0$. This completes the 
proof. 

\qed

 On the other hand,
$\forall h\in V, 
X\in E$, we have $Ad(h)X\in E$. So there is a representation
\[
\rho : V\rightarrow Aut(E),\qquad h\mapsto Ad(h)
\]

Suppose $T'$ is the homogeneous bundle
\[
T'=G\times_V E
\]
whose local section can be represented as $C^\infty$ functions
\[
f: G\rightarrow E
\]
which is $V$ equivariant
\[
f(ga)=Ad(a^{-1})f(g)
\]
for $a\in V, g\in G$. Our next lemma is

\begin{lemma}
\[
T'=T_h(D)
\]
\end{lemma}

{\bf Proof:}
What we are going to prove is that both vector bundles will be coincided 
as subbundles of $T(D)$.

Suppose $\xi\in T'_{gV}$ for $g\in G$ where $T'_{gV}$ is the fiber
of $T'$ at $gV$. Then 
$\xi$ can be represented as 
\[
\xi=(g,\xi_1)\qquad for \quad \xi_1\in E
\]

So the 1-jet in the $\xi$ direction is 
$(g+\varepsilon g\xi_1)V$ for $\varepsilon$ small. Such a point is
\[
(g+\varepsilon g\xi_1)\{f^p\}=(1+\varepsilon g\xi_1 g^{-1})\{F^p\}
\]
where $\{F^p\}=g\{f^p\}$.

Suppose $\xi_2=g\xi_1 g^{-1}$, then
\[
\xi_2 F^p\subset F^{p-1}
\]

Thus $\xi_2\in(T_h)_{gV}(D)$ and
\[
T'_{gV}\subset (T_h)_{gV}(D)
\]

Thus
\[
T'\subset T_h(D)
\]
and $T'$ is the subbundle of $T_h(D)$. But since they 
coincides at the origin, they are equal.

\begin{cor}
Suppose $T_v(D)$ is the distribution of the tangent vectors of the fibers 
of the natural projection
\[
p:\, D\rightarrow G/K
\]
then
\[
T_v(D)\cap T_h(D)=\{ 0\}
\]
\end{cor}

{\bf Proof:}
\[
T_v(D)=G\times_V{\frk v}_1
\]
where ${\frk f}_0={\frk v}+{\frk v}_1$ and ${\frk v}_1$ is the orthonormal 
complement of the Lie algebra ${\frk v}$ of $V$.
 
\qed

\begin{definition}\label{hors}
Let ${U}$ be a complex manifold. If ${U}\subset D$
is a complex submanifold such that $T({U})\subset T_h(D)|_{U}$.
Then we say that ${U}$ is a horizontal slice. If
\[
f:{U}\rightarrow D
\]
is an immersion and $f({U})$ is a horizontal slice, then
we say that $(U,f)$ or $U$ is a horizontal slice. In a word,
a horizontal slice ${U}$ of $D$ is a complex integral submanifold
of the distribution $T_h(D)$.
\end{definition}

Because to become a horizontal slice is a local property, we make the
following definition:

\begin{definition}
Suppose $\Gamma$ is a discrete subgroup of $U$ and suppose $\Gamma\subset
G$ for $D=G/V$. Then if $U\rightarrow D$ is a horizontal slice, we 
also say that $\Gamma\backslash U$ is a horizontal slice.
\end{definition}

\begin{cor}\label{imm}
If $f: U\rightarrow D$ is a horizontal slice, then
\[
p: {U}\subset D\rightarrow G/K_0
\]
is an immersion, where $p: D\rightarrow G/K_0$ is the natural projection
in Definition~\ref{pro}.
\end{cor}

\qed

\section{A Metric Rigidity Theorem}

In this section, we prove that, for concave horizontal slices,
the Hodge metric is intrinsically defined. That is, the Hodge metric does
not
depend on the immersion to the classifying space. 

To be precise, suppose $\Gamma\backslash{U}\rightarrow\Gamma\backslash D$
is a
horizontal slice. Then we can define the Hodge
 metric  on $\Gamma\backslash{U}$. But as
a complex manifold, the horizontal immersion
$\Gamma\backslash{U}\rightarrow\Gamma\backslash D$
may not be unique. If a metric defined on
$\Gamma\backslash{U}$ is independent of the choice of the 
immersion,
we say such a metric is
defined intrinsically.

For the moduli space of a Calabi-Yau threefold, the Hodge metric is
defined intrinsically by the main result in~\cite{Lu5}. It is interesting
to ask if the property is true for  general horizontal slices.

\begin{definition}
The classifying space $D$, as a homogeneous complex manifold, has a
natural invariant K\"ahler form $\omega_H$. In general, $d\omega_H\neq 0$.
However, if $U\rightarrow D$ is a horizontal slice, then $d\omega_H=0$
(cf. ~\cite{Lu6}).
The metric $\omega_H|_U$ is called the Hodge metric.
\end{definition}

\begin{definition}
We say a complex manifold $M$ is {\it concave}, if there is an exhaustion
function $\ph$ on $M$ such that the Hessian of $\ph$ has at least two
negative eigenvalues at each point outside some compact set.
\end{definition}

Any pluriharmonic function on a concave manifold is a constant.

Suppose $f_i: {U}\rightarrow D$, $i=1,2$ are two horizontal
slices. Suppose we have $\Gamma\in Aut({U}), \Gamma_0\in Aut
\, D$ and we have the group homomorphism
\[
\rho : \Gamma\rightarrow \Gamma_0,\qquad 
\]
such that
\[
f_i(\gamma x)=\rho(\gamma)f_i(x),\qquad i=1,2,\quad
\gamma\in\Gamma, x\in U
\]
where the action $\rho(\gamma)$ on $D$ is  the left translation.

The main results of this section are the following two theorems:

\begin{theorem}\label{VHS1}
With the notations as above, suppose that $\Gamma\backslash{U}$ has no
nonconstant pluriharmonic 
functions. Then 
there is an isometry
$f:  f_1(U)\rightarrow  f_2(U)$ such that
 $f\circ  f_1= f_2$.
\end{theorem}

{\bf Proof of Theorem ~\ref{VHS1}:} Let $D_1=G/K_0$ be the  symmetric
space defined in Lemma~\ref{lem34}.
We  denote
$\tilde f_1 :{U}\rightarrow G/K_0$ and $\tilde 
f_2 :{U}\rightarrow G/K_0$ to
be the  two natural projections, that is $\tilde f_i=p\circ f_i$ where
$p$ is defined in Definition~\ref{pro}.
By Corollary~\ref{imm}, both maps are immersions.
 Let
\[
g:{U}\rightarrow \R,\qquad
g(x)=d(f_1(x),f_2(x))
\]
where $d(\cdot,\cdot)$ is the distance function of 
$G/K_0$. Thus since $G/K_0$ is a 
Cartan-Hardamad manifold, $g(x)$ is smooth if $g(x)\neq 0$.

Let $p\in{U}$ and $X\in T_p{U}$. Let 
$X_1=(\tilde f_1)_{*p}X,X_2=(\tilde f_2)_{*p}X$. Let $\sigma$ be the
geodesic ray
starting at $p$ with vector $X$. i.e.
\[
\left\{
\begin{array}{l}
\sigma''(t)=0\\
\sigma(0)=p,\sigma'(0)=X
\end{array}
\right.
\]

Suppose 
the smooth function $\sigma(s,t)$ is defined as follows: for fixed $s$, 
$\sigma(s,t)$ is the geodesic in $G/K$ connecting 
$\tilde f_1(\sigma(s))$ and $\tilde f_2(\sigma(s))$. Furthermore, we
assume 
that $\sigma(0,t)$ is normal. i.e. $t$ is the arc length. define
\[
\tilde X(s)=\frac{d}{ds}{\bigg |}_{s=0} \sigma(s,t)
\]
be the Jacobi field of the variation. In particular
\[
\left\{
\begin{array}{l}
\tilde X(0)=X_1\\
\tilde X(l)=X_2
\end{array}
\right.
\]
where $l=g(x)$. Suppose $T$ is the tangent vector of 
$\sigma(0,t)$, we have the second variation formula
\begin{align*}
XX(g)|_p &=<\nabla_{X_2}X_2, T>-<\nabla_{X_1}X_1, T>\\
& + \int_0^l|\nabla_T\tilde X|^2
-R(T,\tilde X, T, \tilde X)-(T<\tilde X,T>)^2
\end{align*}
where $\nabla$ is the connection operator on $G/K_0$ and 
$R(\cdot,\cdot,\cdot,\cdot)$ is the curvature tensor.

We also have the first variation formula
\[
Xg=<(\tilde f_2)_*X,T>-<(\tilde f_1)_*X,T>
\]

By~\cite[Theorem 1.1]{Lu6}, we know $f_i (i=1,2)$ are pluriharmonic.
That is, we have the following
\[
\nabla_{(\tilde f_i)_*X}(\tilde f_i)_*X+\nabla_{(\tilde
f_i)_*JX}(\tilde f_i)_*JX
+(\tilde f_i)_*J[X,JX]=0
\]
for $i=1,2$.

Define 
\[
D(X,X)=XXg+(JX)(JX)g+J[X,JX]g
\]
Using the fact that $J$ is $\nabla$-parallel, we see

\begin{align*}
D(X,X)g=&\int_0^l|\tilde X'|^2-R(T,\tilde X,T,\tilde X)
-(T<\tilde X,T>)^2\\
+&\int_0^l|\tilde{JX}'|^2-R(T,\widetilde{JX},T,\widetilde{JX})
-(T<\widetilde{JX},T>)^2
\end{align*}
where $\widetilde{JX}$ is the Jacobi connecting $\tilde f_1(J\sigma(t))$ 
and $\tilde f_2(J\sigma(t))$.

{\bf Claim:} If $g(x)\neq 0$, then Hessian of $g$ at $x$ is semipositive.

{\bf Proof:}
Let $(\frac{\partial}{\partial z^1},\cdots, 
\frac{\partial}{\partial z^n})$ be the holomorphic normal 
frame at $p\in{U}$. In order to prove $g$ is 
plurisubharmonic, it suffices to prove that 
$\frac{\partial^2 g}{\partial z_i\partial \bar{z_i}}\geq 0$. But
\[
4\frac{\partial^2\,g}{\partial z_i\partial \bar{z_i}}
=\frac{\partial^2g}{\partial x_i^2}
+\frac{\partial^2 g}{\partial y_i^2}=D(\frac{\partial}{\partial 
x_i},\frac{\partial}{\partial x_i})g
\]

Let $X=\frac{\partial}{\partial x_i}$ in the second variation 
formula. Since the curvature of the symmetric space is nonpositive,
\begin{align*}
&D(\frac{\partial}{\partial x_i},\frac{\partial}{\partial x_i})g
\geq \int_0^l|\tilde X'|^2-(T<\tilde X,T>)^2+|\widetilde{JX}'|^2
-(T<\widetilde{JX},T>)^2
\end{align*}
because the curvature operator is nonpositive. On the other hand
\begin{align}\label{g}
\begin{split}
& |X'|^2-(T<\tilde X,T>)^2=|\tilde X'-<\tilde X',T>T|^2\geq 0\\
& |JX'|^2-(T<\widetilde {JX},T>)^2=|\widetilde {JX}'-<\widetilde 
{JX}',T>T|^2\geq 0
\end{split}
\end{align}
Thus $g$ is plurisubharmonic if $g(x)\neq 0$.

 $g^2(x)$ is a smooth function on ${U}$. It is easy to see
that $g^2$ is a plurisubharmonic function. But $g^2$ is also
$\Gamma$-invariant so it descends to a function on
$\Gamma\backslash{U}$. Thus $g^2$ and $g$ must be 
constant.

Since $g$ is a constant, by Equation~\eqref{g} and the second variational
formula, we
have
\[
\left\{
\begin{array}{l}
\tilde X'-<\tilde X',T>T=0\\
R(T,\tilde X,T,\tilde X)=0
\end{array}
\right.
\]
Moreover, by the first variational formula,
\[
<\tilde X,T>(0)=<\tilde X,T>(l)
\]
Since $\tilde X$ is a Jacobi field, 
  $\tilde
X''\equiv 0$. Furthermore, by the above equations, we have $\tilde
X'\equiv 0$.

This proves that there is an isometry
\[
\tilde f: \tilde f_1(U)\rightarrow\tilde f_2(U),\qquad
\tilde f_1(x)\mapsto \tilde f_2(x)
\]
which sends $\tilde f_1(x)$ to $\tilde f_2(x)$ and thus we have
$\tilde f\circ \tilde f_1=\tilde f_2$.

The theorem follows from the fact that $f_i(U)$ and $\tilde f_i(U)$ are
isometric for $i=1,2$.

\qed

\section{Local Rigidity of the Group Representation}
In this section we study the monodromy group
representation on a horizontal slice. 

We assume that $U$ is a horizontal slice. Let $\Gamma\subset
Aut({U})$ be a discrete group. Suppose $\Gamma\backslash{U}$ is of finite
volume with respect to the Hodge metric.

For the sake of simplicity, we assume that $\Gamma$  is also the subgroup
of the left translation of $D=G/V$, the classifying space. There is a
natural map $\Gamma\backslash{U}\rightarrow\Gamma\backslash G/K_0$
where $G/K_0$ is the  symmetric space
of $D=G/V$  as in Definition~\ref{pro}. Let
\[
{\mathcal G}=\{a\in G | a\in Aut(U)\}
\]

 Let ${\mathcal  G}_0$ be the identity component of ${\mathcal G}$. 

The main theorem of this section is
\begin{theorem}\label{thm51}
Let $\Gamma\backslash{U}$ be of finite Hodge volume. 
Suppose further that ${\mathcal G}_0$
is semisimple 
and ${\mathcal G}_0/{\mathcal K}_0$ is a Hermitian symmetric space but is
not a complex ball, where ${\mathcal K}_0$ is the maximum compact subgroup
of ${\mathcal G}_0$.  Then the representation
 $\Gamma\rightarrow\,{\mathcal G}$ is 
locally rigid.
\end{theorem}

By local rigidity we mean that if $\rho_t:\Gamma\rightarrow {\mathcal G}$
is a
continuous set of representations for $t\in (-\eps,\eps)$,
then there is an $a_t$
for any $|t|<\eps$ such that
$\rho_t=Ad(a_t)\rho_0$.

Before proving the rigidity theorem, we make the following assumption. We
postpone
the proof of the assumption to the end of this section.

\begin{assumption}\label{ass}
Let ${\mathcal K}_0$ be a maximal compact subgroup of ${\mathcal G}_0$.
Suppose
$\Gamma_1=
\Gamma\cap{\mathcal G}_0$. We assume that
$\Gamma_1\backslash{\mathcal G}_0/{\mathcal K}_0$ has finite volume
with respect to the standard Hermitian metric on ${\mathcal
G}_0/{\mathcal K}_0$. In
this case, we will call $\Gamma_1$ has  finite covolume.
\end{assumption}

We prove a series of lemmas.

Let
\[
{\mathcal G}_1=\Gamma +{\mathcal G}_0
\]
be the group generated by $\Gamma$ and ${\mathcal G}_0$ in $G$.

Let
\[
\Gamma_1=\Gamma\cap {\mathcal G}_0
\]

\begin{lemma}\label{lem51}
Let
$\pi : {U}\rightarrow \Gamma\backslash {U}$
be the projection. Then for any $x\in{U}$, the projection of the
${\mathcal 
G}_0$ orbit  $\pi({\mathcal G}_0x)$ is a closed, locally connected, properly embedded 
smooth submanifold of $\Gamma\backslash{U}$.
\end{lemma}

{\bf Proof(cf,~\cite{SF}):}  ${\mathcal G}_0x$ is a closed properly
embedded, locally connected smooth 
submanifold of ${U}$, we claim:

\smallskip

{\bf Claim:} $\pi^{-1}(\pi({\mathcal  G}_0x))={\mathcal G}_1x$.

{\bf Proof:} We know that ${\mathcal G}\subset N({\mathcal G}_0)$, the normalizer of 
${\mathcal G}_0$ in ${\mathcal G}$. So 
$\forall \xi\in {\mathcal G}_0, b\in \Gamma$, there 
is a $\eta\in{\mathcal G}_0$ such that
$b\xi=\eta b$.
Thus $\forall g\in {\mathcal G}_1$, $g=g_1g_2$ where $g_1\in\Gamma$ and $g_2\in{\mathcal 
G}_0$. So
\[
\pi(gx)=\pi(g_1g_2x)=\pi(g_2x)\in\pi({\mathcal G}_0x)
\]

Thus
$gx\in\pi^{-1}\pi({\mathcal G}_0x)$.

On the other hand, if $y\in\pi^{-1}\pi({\mathcal G}_0x)$, 
then $\pi(y)\in\pi({\mathcal 
G}_0x)$, thus by definition, $y\in{\mathcal G}_1x$.

Since ${\mathcal G}_1x$ is a properly embedded, locally connected 
smooth submanifold of ${U}$ and ${\mathcal G}_1x$ is 
$\Gamma$ invariant. The lemma is proved by observing 
$\pi({\mathcal G}_1x)=\pi({\mathcal G}_0x)$.

 \qed

In order to prove Theorem~\ref{thm51}, we use the following famous
theorem of Margulis~\cite{Mar} about the superrigidity of symmetric
spaces:

\begin{thm}[Margulis]\label{margulis}
Suppose that ${\mathcal G}_0$ is defined as above.
If $\Gamma_1$ is of finite covolume, then for any 
homomorphism
\[
\varphi: \Gamma_1\rightarrow\Gamma_1
\]
there is a unique extension
\[
\tilde\varphi : {\mathcal G}_0\rightarrow {\mathcal G}_0
\]
of group homomorphism.
\end{thm}

The following lemma is a straightforward consequence of the above theorem
of Margulis.

\begin{lemma} If $x\in{\mathcal G}_0$ such that
\[
xy=yx
\]
for all $y\in\Gamma_1$, then $x=e$.
\end{lemma}
  
{\bf Proof:}
Let
$\varphi : \Gamma_1\rightarrow\Gamma_1$
by
$y\rightarrow xyx^{-1}$.
Then $\varphi$ has an extension
$\tilde\varphi : {\mathcal G}_0\rightarrow{\mathcal G}_0$.
This extension is unique.
So we must have $\tilde \varphi(y)=xyx^{-1}=y$. Since
${\mathcal G}_0$ is semisimple, we have $x=e$.
\qed

\bigskip

\begin{lemma}
Let $\Gamma_1=\Gamma\cap{\mathcal G}_0$, then 
\[
Out(\Gamma_1)/Inn(\Gamma_1)
\]
is a finite group.
\end{lemma}

Here $Out(\Gamma_1)$ denotes the group of isomorphisms of $\Gamma_1$ and
$Inn(\Gamma_1)$ denotes the group of conjugations of $\Gamma_1$.

{\bf Proof:}
Let 
$\varphi: \Gamma_1\rightarrow\Gamma_1$
be an element in $Out(\Gamma_1)$. 
Since $\Gamma_1$ has finite covolume, we know there is a unique extension
$\tilde\varphi : {\mathcal G}_0\rightarrow {\mathcal G}_0$.

Thus $\tilde\varphi\in Out({\mathcal G}_0)$.  Since $Out({\mathcal 
G}_0)=Inn({\mathcal G}_0)$ 
because ${\mathcal G}_0$ is semisimple, there is a $b\in{\mathcal 
G}_0$ such that $\tilde\varphi(x)=bxb^{-1}$. 
Define 
\[
\tilde\ph: {\mathcal G}_0/{\mathcal K}_0
\rightarrow
{\mathcal G}_0/{\mathcal K}_0
\qquad
a{\mathcal K}_0\rightarrow ba{\mathcal K}_0
\]
It is a $\Gamma$-equivariant holomorphic map. The lemma 
then follows from
the following proposition.

\begin{prop}
Suppose $\Gamma\backslash G/K$ is of finite volume, then
$Aut(\Gamma\backslash G/K)$ is a finite group.
\end{prop}

{\bf Proof:} Since $\Gamma\backslash G/K$ is a Hermitian symmetric space,
we know
$Aut(\Gamma\backslash G/K)$ is the same as
$Iso(\Gamma\backslash G/K)$.

Suppose $Iso(\Gamma\backslash G/K)$ is not finite. Then we have
a sequence of isometries $f_1,f_2,\cdots$. Let
$p\in\Gamma\backslash G/K$ be a fixed point and let $V$ be a
normal coordinate neighborhood of $p$. The we know that $\{f_i(p)\}$ must
be
bounded, otherwise there is a subsequence of $f_i$ such that
$f_i(U)$ will be mutually disjoint. This will contradict to the
fact that $\Gamma\backslash G/K$ has finite volume, because
\[
vol(\Gamma\backslash G/K)\geq\sum vol(f_i(U))=+\infty
\]
A contradiction.
Let $q=lim\,f_i(p)$. For any $x\in\Gamma\backslash G/K$,
if $i$ is large enough such that $d(f_i(p),q)<1$, then
\[
d(f_i(x),q)\leq d(f_i(x),f_i(p))+1=d(x,p)+1
\]
By Ascoli theorem, there is a subsequence of $f_i$ such that
$f_i$ converges to an $f\in Iso(\Gamma\backslash G/K)$.
Thus $Iso(\Gamma\backslash G/K)$ is not discrete.
So there is a holomorphic vector field $X$ on $\Gamma\backslash G/K$.

Suppose $X=X^i\frac{\pa}{\pa z^i}$ in local coordinate,
and  $||X||^2=G_{i\bar j}X^i\bar{X^j}$. Suppose
the local coordinate is normal, then
\begin{equation}\label{X1}
\pa_k\bar{\pa}_l||X||^2
=R_{i\bar jk\bar l} X^i\bar{X}^j
+\pa_kX^i\bar{\pa_lX^i}
\end{equation}
where $R_{i\bar jk\bar l}$ is the curvature tensor of the 
symmetric space. Thus in particular $\pa\bar\pa||X||^2\geq 0$.

By the  theorem of ~\cite{Bor}, $\Gamma\backslash G/K$ is a concave
manifold.
Thus $||X||^2$ is a constant. On the other hand,
from equation  ~\eqref{X1}, we have
\[
\Delta||X||^2=Ric(X)+|\nabla X|^2
\]

So $Ric(X)=0$ and thus $X\equiv 0$. This is a contradiction.

\qed

From the above proposition,
there is an integer $n$ such that  $a^n=e$ for all
\[
a\in Out(\Gamma_1)/Inn(\Gamma_1)
\]

Let $\tilde\Gamma$ be the subgroup 
of ${\mathcal G}_1$ generated by $\Gamma_1$ and $a^n$ 
where
$a\in\Gamma$. 
Then we have an exact sequence
\begin{equation}\label{ext}
1\rightarrow\Gamma_1\rightarrow \tilde\Gamma\rightarrow\tilde B\rightarrow 1
\end{equation}
where $\tilde B$ is the quotient $\tilde \Gamma/\Gamma_1$. For any 
$\bar{b}\in\tilde\Gamma/\Gamma_1$ with $b\in\tilde\Gamma$, we have
$b\Gamma_1b^{-1}\subset\Gamma_1$. So $b\in Out(\Gamma_1)$. But by the definition of 
$\tilde\Gamma$, $b$ is a trivial element in $Out(\Gamma_1)/Inn(\Gamma_1)$. So there is a 
$c\in\Gamma_1$ such that $bc$ is commutative to $\Gamma_1$. So 
there is a homomorphism
\[
\eta : B\rightarrow \tilde\Gamma,\qquad \bar{b}\mapsto bc
\]
We can thus define a homomorphism
\[
\xi : \Gamma_1\times\tilde \Gamma/\Gamma_1\rightarrow\tilde \Gamma
\]
such that
\[
\xi(a,b)=a\eta(b)
\]
which is an isomorphism. In other words,
the exact sequence~\eqref{ext} splits.

\begin{lemma}
Let $\tilde{\mathcal G}_1=\tilde\Gamma+{\mathcal G}_0$. Then
\[
\tilde{\mathcal G}_1={\mathcal G}_0\times \tilde B
\]
\end{lemma}

{\bf Proof:}
We have
\[
\tilde\Gamma=\Gamma_1\times \tilde B
\]
Define
\[
\varphi :{\mathcal G}_0\times \tilde B\rightarrow 
\tilde{\mathcal G}_1,\quad
\varphi(a,b)=ab
\]
Then $\varphi$ is an isomorphism.

\qed

Thus we know a family of representation of $\tilde\Gamma$ 
splits to the representation 
to the discrete group $\tilde B$ and Lie group 
${\mathcal G}_0$ respectively.  

\begin{lemma}
If the representation
 $\tilde\Gamma\rightarrow\tilde{\mathcal G}_1$ is locally 
rigid, then the representation
$\Gamma\rightarrow\tilde{\mathcal G}$ is
also locally rigid.
\end{lemma}

{\bf Proof:} Let $\ph_t:\Gamma\rightarrow\tilde{\mathcal G}$ 
be a local family of representations, $t\in(-\eps,\eps)$.
Then we see that $\ph_t$ restricts to a trivial family
of representations on $\tilde\Gamma$. That is, there are
$a_t\in\tilde{\mathcal G}_1\subset\tilde{\mathcal G}$
with $a_0=e$ such that $\ph_t(x)=a_t\ph_0(x)a_t^{-1}$
for $x\in\tilde\Gamma$. Let $\xi_t=Ad(a_t^{-1})\ph_t$. 
Then we know $\xi_t(x)=\ph_0(x)$ for all $x\in\tilde\Gamma$.
Now if $x\in\Gamma$, then $x^n\in\tilde\Gamma$.
So we have $(\xi_t(x))^n=(\ph_0(x))^n$
and $\xi_0(x)=\ph_0(x)$. Thus
$\xi_t(x)=\ph_0(x)$.

\qed

In the rest of this section, we prove Assumption~\ref{ass}.

 \begin{lemma}
${\mathcal G}_1$ is a closed subgroup of $G$.
\end{lemma}

{\bf Proof:}
We know that ${\mathcal G}_1\subset{\mathcal G}$. Let $x_m\in{\mathcal G}_1$ such that 
$x_m\rightarrow x$ for $x\in G$. Then $x\in Aut({\mathcal M})$ so
$x\in{\mathcal G}$. 
Thus for sufficient large $m$, $x_m$ and $x$ are in the same component. In particular, 
we have $x\in{\mathcal G}_1$.
\qed

\begin{lemma}
Let $p\in {U}$, we have
\[
\underset{q\in{\mathcal G}_1\backslash{\mathcal G}_0}{inf} 
d(qp,{\mathcal G}_0p)>0
\]
\end{lemma}

{\bf Proof:} Suppose the assertion is not true, then we have $\{q_m\}\in{\mathcal G}_1$
and $g_m\in{\mathcal G}_0$ 
such that
\[
d(q_mp,g_mp)\rightarrow 0,\quad m\rightarrow +\infty
\]
or
\[
d(g_m^{-1}q_mp,p)\rightarrow 0,\quad m\rightarrow +\infty
\]

It is easy to check that ${\mathcal G}_0p$ is a homogeneous manifold, with compact 
stable group. 

Thus, there are $k_m\in{\mathcal K}_0$, a compact subgroup of ${\mathcal G}_0$ such that
\[
g_m^{-1}q_mk_m\rightarrow e
\]
So by passing a subsequence if necessary, we know
\[
g_m^{-1}q_m\rightarrow g\in {\mathcal K}_0\subset{\mathcal G}_0
\]

This contradicts the fact that ${\mathcal G}_0$ is open.

Let $x,y\in U$. Let
\[
L_1={\mathcal G}_0x,\qquad L_2={\mathcal G}_0y
\]
be the two ${\mathcal G}_0$ orbits. We can define
\[
f(p)=d(p,L_2)
\]
be the distance of a point $p\in L_1$ to $L_2$.

\begin{lemma}
$f(p)$ is a constant.
\end{lemma}

{\bf Proof:}
Let $p,q\in L_1$. Then there is a $g\in{\mathcal G}_0$ such that $q=gp$.
We have
\[
d(q,L_2)\leq d(q,\xi)=d(gp,\xi)=d(p,g^{-1}\xi)
\]
This proves
\[
d(q,L_2)\leq d(p,L_2)
\]
On the other hand, we also have
\[
d(q,L_2)\leq d(p,L_2)
\]
Thus $d(q,L_2)=d(p,L_2)$ and $f(p)$ is a constant.

Define the distance between two orbits $L_1,L_2$ by
$d(L_1,L_2)=f(p)$. If $L_1\neq L_2$, $d(L_1,L_2)>0$.

Let $\bar{a}\in{\mathcal G}_1/{\mathcal G}_0$. Then $\bar{a}L$ defines 
another orbit. So 
we have a map
\[
{\mathcal G}_1/{\mathcal G}_0\rightarrow \R,\qquad \bar{a}\mapsto
d(\bar{a}L,L)
\]
We know that $d(\bar aL,L)>0$ for $\bar a\neq 0$. Furthermore we have

\begin{lemma}
\[
\varepsilon=\underset{\bar a\neq 0}{inf}\,d(\bar aL,L)>0
\]
\end{lemma}

{\bf Proof:} This is a consequence of the previous two lemmas.

\qed

For any orbit ${\mathcal G}_0p$, $\Gamma\backslash{\mathcal G}_0p$ is
a closed, 
properly embedded submanifolds
\newline
(Lemma~\ref{lem51}). We fix one of them, say
$L$.

Let 
\[
W=\{x\in{U}|d(x,L)<\frac{\varepsilon}{100}\}
\]
where $\varepsilon$ is defined in the previous lemma. 
Then for any $a\in{\mathcal 
G}_1\backslash{\mathcal G}_0$, $aU\cap U=\emptyset$. In particular
\[
\Gamma\backslash{W}=\Gamma_1\backslash{W}
\]
in $\Gamma\backslash{U}$.

Now that
\[
vol(\Gamma\backslash{U})\geq vol(\Gamma\backslash W)=vol
(\Gamma_1\backslash W)
\]
For any $p\in\Gamma_1\backslash U$, there is a unique $q\in L$ such that
\[
d(p,q)=d(p,L)
\]

Now we can prove
the following proposition which implies the assumption:

\begin{prop}
If $vol(\Gamma\backslash{U})<+\infty$, then
\[
vol(\Gamma\backslash L)<+\infty
\]
\end{prop}

{\bf Proof:}
Let
$f(p)=d(p,L)$. Then by the coarea formula
\[
vol(\Gamma_1\backslash 
U)=\int^{\varepsilon}_0(\int_{f=c}\frac{1}{|\nabla f|})dc
\]
But $|\nabla f|\leq 1$. So
\[
vol(\Gamma_1\backslash U)\geq\int_0^{\epsilon} vol(f=c)dc
\]
so at least there is a $c$ s.t.
\[
vol(f=c)<\infty
\]
Note that $\dim\,\{f=c\}=\dim U-1$. The proposition then follows from
the induction.

\bibliographystyle{acm}
\bibliography{bib}

\end{document}